\documentclass[12pt]{article}       
\usepackage{amsfonts,amsmath,amssymb}
\usepackage[usenames,dvipsnames]{pstricks}
\usepackage{epsfig}
\usepackage{pst-grad}												 
\usepackage{pst-plot} 												

\usepackage{graphicx}

\newtheorem{thm}{Theorem}
\newtheorem{cor}[thm]{Corollary}
\newtheorem{lem}[thm]{Lemma}

\newenvironment{pr}
{\noindent\textbf{Proof.}\ \ }{\hfill $\Box$ \medskip}

\title{Some conditions on 5-cycles that make planar graphs 4-choosable }
\author {Pongpat Sittitrai\\ 
{\small\em Department of Mathematics, Faculty of Science, Khon Kaen University, 40002, Thailand }\\  
{\small\em E-mail address: pongpat@kkumail.com} 
\and Kittikorn Nakprasit \footnote{Corresponding Author} \\ 
{\small\em Department of Mathematics, Faculty of Science, Khon Kaen University, 40002, Thailand }\\
{\small\em E-mail address: kitnak@hotmail.com}}
\date{}

\begin{document}

\maketitle

\begin{center}{\bf Abstract}\end{center}
\indent\indent 
Consider two conditions on a graph:  
(1) each 5-cycle is not a subgraph of 5-wheel 
and does not share exactly one edge with 3-cycle, 
and (2) each 5-cycle is not adjacent to two 3-cycles 
and is not adjacent to a 4-cycle with chord.  
We show that if a planar graph $G$ satisfies one 
of the these conditions, then $G$ is 4-choosable. 
This yields that 
if each 5-cycle of a planar graph $G$ is not adjacent a 3-cycle, 
then $G$ is 4-choosable. 

\section{Introduction}

The concept of list coloring was introduced by Vizing \cite{Vizing} 
and by Erd\H os, Rubin, and Taylor \cite{Erdos}, independently.  
An \emph{assignment} $L$ for a graph $G$ assigns a list $L(v)$ 
for each vertex $v.$ If $|L(v)|= k$ for each vertex $v,$ then we 
call $L$ a $k$-assignment. 
A graph $G$ is $L$-colorable if we can color $G$ with each vertex has a color 
from its list and no two adjacent vertices receive the same color.   
If $G$ is $L$-colorable for any $k$-assignment $L,$ then we say  $G$ is $k$-\emph{choosable}. 

Thomassen \cite{Tho} proved that every planar graph is 5-choosable. 
Voight \cite{Vo1} and Mirzakhani \cite{Mir} 
presented the examples of non 4-choosable graphs.
Additionally, Gutner \cite{Gut} showed that the problem of determining whether 
a planar graph is 4-choosable is NP-hard.  
This leads to the interest of finding some nice sufficient conditions 
for planar graphs to be 4-choosable. 
If a planar graph has no 3-cycles, then it has a vertex of degree at most 3. 
Thus it is 4-choosable. 
More sufficient conditions are founded, 
for examples, it was shown that a planar graph is 4-choosable if it has no 
4-cycles \cite{Lam2}, 5-cycles \cite{Lam1, Wang1}, 6-cycles \cite{Fi}, 
7-cycles \cite{Far}, intersecting 3-cycles \cite{Wang2}, 
intersecting 5-cycles \cite{Hu},  
3-cycles adjacent to 4-cycles \cite{Bo, Cheng},  
5-cycles simultaneously adjacent to 3-cycles and 4-cycles. 

In this paper, we prove the following theorems. 

\begin{thm}\label{main1} 
If each 5-cycle of a planar graph $G$ is not a subgraph of 5-wheel 
and does not share exactly one edge with 3-cycle, 
then $G$ is 4-choosable. 
\end{thm}

\begin{thm}\label{main2} 
If each 5-cycle of a planar graph $G$ is not adjacent to two 3-cycles 
and is not adjacent to a 4-cycle with chord, then $G$ is 4-choosable. 
\end{thm} 

The following corollary is an easy consequence. 
\begin{cor}
If each 5-cycle of a planar graph $G$ is not adjacent a 3-cycle, 
then $G$ is 4-choosable. 
\end{cor} 


\section{Preliminaries} 

First, we introduce some notations and definitions. 
A $k$-vertex (face) is a vertex (face) of degree $k,$ 
a $k^+$-vertex (face) is a vertex (face) of degree at least $k,$ 
and a $k^-$-vertex (face) is a vertex (face) of degree at most $k.$ 
A \emph{$(d_1,d_2,\dots,d_k)$-face} $f$ is a face of degree $k$ 
where all vertices on $f$ have degree $d_1,d_2,\dots,d_k$. 
A \emph{$(d_1,d_2,\dots,d_k)$-vertex} $v$ is a vertex of degree $k$ 
where all faces incident to $v$ have degree $d_1,d_2,\dots,d_k$. 

A \emph{trio} is a graphs consist of a vertex set of five elements, 
namely $\{x,y,z, u, v, w\},$ and an edge set $\{xy, xu, xv, yv, yw, uv, vw\}$ 
(see Fig. \ref{fig1}). 
We call $s$ on a 3-cycle $f$ a \emph{good} vertex of $f$ if $f$ is not in any trio, 
a \emph{bad} vertex of $f$ if $f$ is in a trio and is the only 3-cycle 
on that trio that $s$ is belong to, 
a \emph{worst} vertex of $f$ if $s$ is on $f$ and all other 3-cycles, 
otherwise $s$ is called a \emph{worse} vertex of $f.$  
We call a face $f$ is a \emph{bad} (\emph{worse}, or \emph{worst}, respectively) face  
of a vertex $v$ if $v$ is a bad (worse, or worst, respectively) vertex of  $f.$

\begin{figure}[ht]\label{fig1}
\centering
\scalebox{1} 
{
\begin{pspicture}(0,-1.5829687)(5.2228127,1.5829687)
\psdots[dotsize=0.2](1.5809375,0.88453126)
\psdots[dotsize=0.2](2.5809374,-0.91546875)
\psdots[dotsize=0.2](3.5809374,0.88453126)
\psdots[dotsize=0.2](4.5809374,-0.91546875)
\psdots[dotsize=0.2](0.5809375,-0.91546875)
\psline[linewidth=0.04cm](0.5809375,-0.91546875)(1.5809375,0.88453126)
\psline[linewidth=0.04cm](1.5809375,0.88453126)(3.5809374,0.88453126)
\psline[linewidth=0.04cm](0.5809375,-0.91546875)(4.5809374,-0.91546875)
\psline[linewidth=0.04cm](4.5809374,-0.91546875)(3.5809374,0.88453126)
\psline[linewidth=0.04cm](1.5809375,0.88453126)(2.5809374,-0.91546875)
\psline[linewidth=0.04cm](2.5809374,-0.91546875)(3.5809374,0.88453126)
\usefont{T1}{ptm}{m}{n}
\rput(1.4223437,1.3945312){$x$}
\usefont{T1}{ptm}{m}{n}
\rput(3.8323438,1.3945312){$y$}
\usefont{T1}{ptm}{m}{n}
\rput(0.23234375,-1.4054687){$u$}
\usefont{T1}{ptm}{m}{n}
\rput(2.6323438,-1.4054687){$v$}
\usefont{T1}{ptm}{m}{n}
\rput(4.8723435,-1.4054687){$w$}
\end{pspicture} 
}
\caption{A trio}
\end{figure}

\begin{figure}[ht]\label{fig2}
\centering
\scalebox{1} 
{
\begin{pspicture}(0,-2.4076562)(16.542812,2.4076562)
\psdots[dotsize=0.2](1.9809375,1.7092187)
\psdots[dotsize=0.2](2.9809375,-0.09078125)
\psdots[dotsize=0.2](3.9809375,1.7092187)
\psdots[dotsize=0.2](4.9809375,-0.09078125)
\psdots[dotsize=0.2](0.9809375,-0.09078125)
\psline[linewidth=0.04cm](0.9809375,-0.09078125)(1.9809375,1.7092187)
\psline[linewidth=0.04cm](1.9809375,1.7092187)(3.9809375,1.7092187)
\psline[linewidth=0.04cm](0.9809375,-0.09078125)(4.9809375,-0.09078125)
\psline[linewidth=0.04cm](4.9809375,-0.09078125)(3.9809375,1.7092187)
\psline[linewidth=0.04cm](1.9809375,1.7092187)(2.9809375,-0.09078125)
\psline[linewidth=0.04cm](2.9809375,-0.09078125)(3.9809375,1.7092187)
\usefont{T1}{ptm}{m}{n}
\rput(2.0323439,2.2192187){$d(x)\leq 5$}
\usefont{T1}{ptm}{m}{n}
\rput(4.2323437,2.2192187){$d(y)=4$}
\usefont{T1}{ptm}{m}{n}
\rput(0.63234377,-0.5807812){$d(u)=4$}
\usefont{T1}{ptm}{m}{n}
\rput(3.0323439,-0.5807812){$d(v)=4$}
\usefont{T1}{ptm}{m}{n}
\rput(5.2723436,-0.5807812){$d(w)=4$}
\psdots[dotsize=0.2](8.380938,1.7092187)
\psdots[dotsize=0.2](9.380938,-0.09078125)
\psdots[dotsize=0.2](10.380938,1.7092187)
\psdots[dotsize=0.2](11.380938,-0.09078125)
\psdots[dotsize=0.2](7.3809376,-0.09078125)
\psline[linewidth=0.04cm](7.3809376,-0.09078125)(8.380938,1.7092187)
\psline[linewidth=0.04cm](8.380938,1.7092187)(10.380938,1.7092187)
\psline[linewidth=0.04cm](7.3809376,-0.09078125)(11.380938,-0.09078125)
\psline[linewidth=0.04cm](11.380938,-0.09078125)(10.380938,1.7092187)
\psline[linewidth=0.04cm](8.380938,1.7092187)(9.380938,-0.09078125)
\psline[linewidth=0.04cm](9.380938,-0.09078125)(10.380938,1.7092187)
\usefont{T1}{ptm}{m}{n}
\rput(8.222343,2.2192187){$x$}
\usefont{T1}{ptm}{m}{n}
\rput(10.632343,2.2192187){$y$}
\usefont{T1}{ptm}{m}{n}
\rput(7.032344,-0.5807812){$u$}
\usefont{T1}{ptm}{m}{n}
\rput(9.4323435,-0.5807812){$v$}
\usefont{T1}{ptm}{m}{n}
\rput(11.672344,-0.5807812){$w$}
\psdots[dotsize=0.2](13.980938,-0.09078125)
\psdots[dotsize=0.2](15.980938,-0.09078125)
\psdots[dotsize=0.2](14.980938,1.7092187)
\psline[linewidth=0.04cm](13.980938,-0.09078125)(15.980938,-0.09078125)
\psline[linewidth=0.04cm](13.980938,-0.09078125)(14.980938,1.7092187)
\psline[linewidth=0.04cm](14.980938,1.7092187)(15.980938,-0.09078125)
\usefont{T1}{ptm}{m}{n}
\rput(15.232344,2.0192187){$y$}
\usefont{T1}{ptm}{m}{n}
\rput(13.832344,-0.5807812){$v$}
\usefont{T1}{ptm}{m}{n}
\rput(16.072344,-0.5807812){$w$}
\usefont{T1}{ptm}{m}{n}
\rput(8.052343,1.8592187){$2$}
\usefont{T1}{ptm}{m}{n}
\rput(7.072344,0.07921875){$2$}
\usefont{T1}{ptm}{m}{n}
\rput(11.632343,0.09921875){$2$}
\usefont{T1}{ptm}{m}{n}
\rput(10.652344,1.8792187){$3$}
\usefont{T1}{ptm}{m}{n}
\rput(9.372344,0.33921874){$4$}
\usefont{T1}{ptm}{m}{n}
\rput(14.632343,2.0192187){$2$}
\usefont{T1}{ptm}{m}{n}
\rput(13.632343,0.21921875){$2$}
\usefont{T1}{ptm}{m}{n}
\rput(16.232344,0.21921875){$2$}
\usefont{T1}{ptm}{m}{n}
\rput(2.9282813,-1.3557812){\large $H$}
\usefont{T1}{ptm}{m}{n}
\rput(2.9523437,-1.9807812){$(1)$}
\usefont{T1}{ptm}{m}{n}
\rput(9.352344,-1.9807812){$(2)$}
\usefont{T1}{ptm}{m}{n}
\rput(14.952344,-1.9807812){$(3)$}
\end{pspicture} 
}
\caption{Subgraphs of $G$ and their related list assignments}
\end{figure}
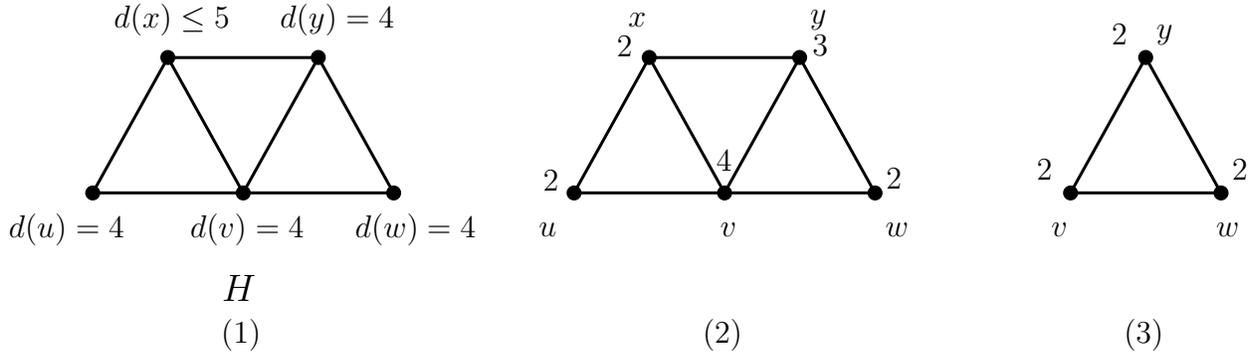

\begin{lem}\label{lem1}
Suppose $G$ is a minimal non 4-choosable graph. 
Then each of the followings holds.\\  
(1) Each vertex degree is at least 4. \\
(2) $G$ does not contain a $(4,4,4,4)$-face. \\
(3) $G$ does not contain a subgraph $H$ as in Fig. \ref{fig2}(1). 
\end{lem} 

\begin{pr} Let $L$ be a 4-list assignment of $G.$\\ 
(1) This is a well-known fact, see e.g. \cite{Lam2}.\\ 
(2) Suppose $G$ contains a $(4,4,4,4)$-face $F.$ 
By minimality of $G,$ the graph $G -V(F)$ has an $L$-coloring. 
A list of available colors of each vertex in $F$ 
has size at least 2. 
Thus we can extend $L$-coloring to $G,$ a contradiction.\\ 
(3) Suppose $G$ contains  a subgraph $H$ as in Fig. \ref{fig2}. 
By minimality of $G,$ the graph $G -V(H)$ has an $L$-coloring. 
The lower bound on size of each residual list is as in Fig. \ref{fig2}(2). 
Next, choose colors for $x$ and $u.$ 
We can extend an $L$-coloring to $G$ unless every residual list is 
identical and has size 2 but one can recolor $x$ and $u$ to avoid this situation. 
This completes the proof.   
\end{pr} 

Next, we introduce a famous theorem established by Alon and Tarsi \cite{Alon}. 
This theorem shows the connection between 
the list coloring and an orientation of a graph. 
A digraph $D$ consists of a vertex set $V(D)$ and an arc set $A(D).$ 
For an arc $a =(u,v) \in A(D),$ we say $a$ has a tail $u$ and a head $v.$ 
The indegree $d^-_D(v)$ of a vertex $v$ is the number of arcs with head $v,$ 
and the outdegree $d^+_D(v)$ of a vertex $v$ is the number of arcs with tail $v.$ 
A subgraph $H$ of $D$ is called \emph{Eulerian} if $d^-_H(v)=d^+_H(v)$ 
for each vertex $v$ of $H.$ 
In this context, $H$ can be disconnected or even edgeless. 
A digraph $H$ is \emph{even} (\emph{odd}) if it has even (odd) number of arcs. 
Let $EE(D)$ denote the number of even spanning Eulerian subgraphs of $D,$ 
and $EO(D)$ denote the number of odd spanning Eulerian subgraphs of $D.$ 

\begin{thm}\label{null}\cite{Alon}
Let $D$ be a digraph. Let $L$ be a list assignment with $|L(v)| =d^+_D(v)+1$ 
for each vertex $v$ of $D.$ If $EE(D) \neq EO(D),$ then $D$ is $L$-colorable. 
\end{thm}

\begin{lem}\label{lem2}
A minimal non 4-choosable graph $G$ does not contain 
a subgraph isomorphic to one of the configurations in Fig. \ref{fig3}. 
\end{lem} 

\begin{figure}[ht]\label{fig3}
\centering
\scalebox{0.9} 
{
\begin{pspicture}(0,-2.1784375)(17.22,2.1584375)
\psdots[dotsize=0.2](1.6,1.5384375)
\psdots[dotsize=0.2](2.6,-0.2615625)
\psdots[dotsize=0.2](3.6,1.5384375)
\psdots[dotsize=0.2](4.6,-0.2615625)
\psdots[dotsize=0.2](0.6,-0.2615625)
\psline[linewidth=0.04cm](0.6,-0.2615625)(1.6,1.5384375)
\psline[linewidth=0.04cm](1.6,1.5384375)(3.6,1.5384375)
\psline[linewidth=0.04cm](0.6,-0.2615625)(4.6,-0.2615625)
\psline[linewidth=0.04cm](4.6,-0.2615625)(3.6,1.5384375)
\psline[linewidth=0.04cm](1.6,1.5384375)(2.6,-0.2615625)
\psline[linewidth=0.04cm](2.6,-0.2615625)(3.6,1.5384375)
\usefont{T1}{ptm}{m}{n}
\rput(2.3714063,-1.9515625){$(1)$}
\usefont{T1}{ptm}{m}{n}
\rput(8.771406,-1.7515625){$(2)$}
\usefont{T1}{ptm}{m}{n}
\rput(14.371407,-1.7515625){$(3)$}
\psdots[dotsize=0.2](6.6,-0.2615625)
\psdots[dotsize=0.2](6.6,1.5384375)
\psdots[dotsize=0.2](8.6,1.5384375)
\psdots[dotsize=0.2](10.6,1.5384375)
\psdots[dotsize=0.2](8.6,-0.2615625)
\psdots[dotsize=0.2](10.6,-0.2615625)
\psdots[dotsize=0.2](12.6,-0.2615625)
\psdots[dotsize=0.2](12.6,1.5384375)
\psdots[dotsize=0.2](14.6,1.5384375)
\psdots[dotsize=0.2](14.6,-0.2615625)
\psdots[dotsize=0.2](16.6,-0.2615625)
\psline[linewidth=0.04cm](6.6,-0.2615625)(10.6,-0.2615625)
\psline[linewidth=0.04cm](12.6,-0.2615625)(16.6,-0.2615625)
\psline[linewidth=0.04cm](6.6,-0.2615625)(6.6,1.5384375)
\psline[linewidth=0.04cm](6.6,1.5384375)(10.6,1.5384375)
\psline[linewidth=0.04cm](10.6,1.5384375)(10.6,-0.2615625)
\psline[linewidth=0.04cm](8.6,-0.2615625)(8.6,1.5384375)
\psline[linewidth=0.04cm](12.6,1.5384375)(12.6,-0.2615625)
\psline[linewidth=0.04cm](12.6,1.5384375)(14.6,1.5384375)
\psline[linewidth=0.04cm](14.6,1.5384375)(14.6,-0.2615625)
\psline[linewidth=0.04cm](14.6,1.5384375)(16.6,-0.2615625)
\psline[linewidth=0.04cm](6.6,1.5384375)(6.6,2.1384375)
\psline[linewidth=0.04cm](6.6,1.5384375)(6.0,1.5384375)
\psline[linewidth=0.04cm](8.6,1.5384375)(8.6,2.1384375)
\psline[linewidth=0.04cm](8.6,-0.2615625)(8.6,-0.8615625)
\psline[linewidth=0.04cm](10.6,1.5384375)(10.6,2.1384375)
\psline[linewidth=0.04cm](10.6,1.5384375)(11.2,1.5384375)
\psline[linewidth=0.04cm](10.6,-0.2615625)(11.2,-0.2615625)
\psline[linewidth=0.04cm](10.6,-0.2615625)(10.6,-0.8615625)
\psline[linewidth=0.04cm](6.6,-0.2615625)(6.0,-0.2615625)
\psline[linewidth=0.04cm](6.6,-0.2615625)(6.6,-0.8615625)
\psline[linewidth=0.04cm](1.6,1.5384375)(1.2,1.9384375)
\psline[linewidth=0.04cm](3.6,1.5384375)(4.0,1.9384375)
\psline[linewidth=0.04cm](4.6,-0.2615625)(5.2,-0.2615625)
\psline[linewidth=0.04cm](4.6,-0.2615625)(4.6,-0.8615625)
\psline[linewidth=0.04cm](0.6,-0.2615625)(0.0,-0.2615625)
\psline[linewidth=0.04cm](0.6,-0.2615625)(0.6,-0.8615625)
\psline[linewidth=0.04cm](0.6,-0.2615625)(0.2,-0.6615625)
\psline[linewidth=0.04cm](12.6,1.5384375)(12.6,2.1384375)
\psline[linewidth=0.04cm](12.6,1.5384375)(12.0,1.5384375)
\psline[linewidth=0.04cm](12.6,-0.2615625)(12.0,-0.2615625)
\psline[linewidth=0.04cm](12.6,-0.2615625)(12.6,-0.8615625)
\psline[linewidth=0.04cm](14.6,1.5384375)(14.6,2.1384375)
\psline[linewidth=0.04cm](14.6,-0.2615625)(14.6,-0.8615625)
\psline[linewidth=0.04cm](14.6,1.5384375)(15.2,1.7384375)
\psline[linewidth=0.04cm](16.6,-0.2615625)(17.2,-0.0615625)
\psline[linewidth=0.04cm](16.6,-0.2615625)(17.0,-0.6615625)
\end{pspicture} 
}
\caption{Three Subgraphs of $G$ }
\end{figure}
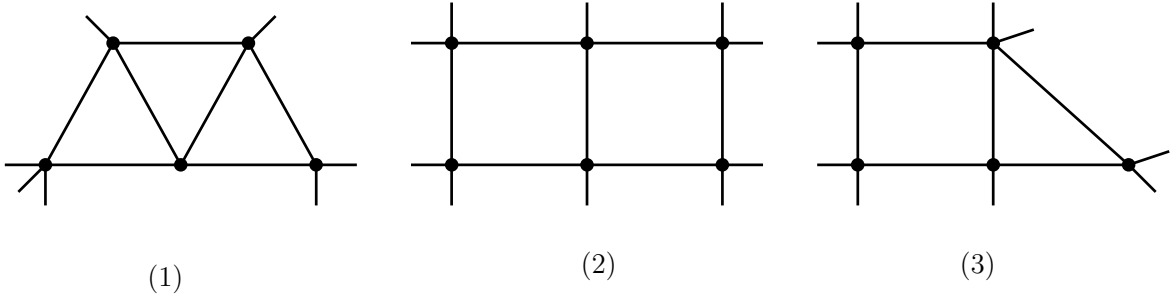
\begin{figure}[ht]\label{fig4}
\centering
\scalebox{0.9} 
{
\begin{pspicture}(0,-2.2067187)(16.782812,2.2067187)
\psdots[dotsize=0.2](1.3809375,1.6082813)
\psdots[dotsize=0.2](2.3809376,-0.19171876)
\psdots[dotsize=0.2](3.3809376,1.6082813)
\psdots[dotsize=0.2](4.3809376,-0.19171876)
\psdots[dotsize=0.2](0.3809375,-0.19171876)
\psline[linewidth=0.04cm](0.3809375,-0.19171876)(1.3809375,1.6082813)
\psline[linewidth=0.04cm](1.3809375,1.6082813)(3.3809376,1.6082813)
\psline[linewidth=0.04cm](0.3809375,-0.19171876)(4.3809376,-0.19171876)
\psline[linewidth=0.04cm](4.3809376,-0.19171876)(3.3809376,1.6082813)
\psline[linewidth=0.04cm](1.3809375,1.6082813)(2.3809376,-0.19171876)
\psline[linewidth=0.04cm](2.3809376,-0.19171876)(3.3809376,1.6082813)
\psdots[dotsize=0.2](6.3809376,-0.19171876)
\psdots[dotsize=0.2](6.3809376,1.6082813)
\psdots[dotsize=0.2](8.380938,1.6082813)
\psdots[dotsize=0.2](10.380938,1.6082813)
\psdots[dotsize=0.2](8.380938,-0.19171876)
\psdots[dotsize=0.2](10.380938,-0.19171876)
\psdots[dotsize=0.2](12.380938,-0.19171876)
\psdots[dotsize=0.2](12.380938,1.6082813)
\psdots[dotsize=0.2](14.380938,1.6082813)
\psdots[dotsize=0.2](14.380938,-0.19171876)
\psdots[dotsize=0.2](16.380938,-0.19171876)
\psline[linewidth=0.04cm](6.3809376,-0.19171876)(10.380938,-0.19171876)
\psline[linewidth=0.04cm](12.380938,-0.19171876)(16.380938,-0.19171876)
\psline[linewidth=0.04cm](6.3809376,-0.19171876)(6.3809376,1.6082813)
\psline[linewidth=0.04cm](6.3809376,1.6082813)(10.380938,1.6082813)
\psline[linewidth=0.04cm](10.380938,1.6082813)(10.380938,-0.19171876)
\psline[linewidth=0.04cm](8.380938,-0.19171876)(8.380938,1.6082813)
\psline[linewidth=0.04cm](12.380938,1.6082813)(12.380938,-0.19171876)
\psline[linewidth=0.04cm](12.380938,1.6082813)(14.380938,1.6082813)
\psline[linewidth=0.04cm](14.380938,1.6082813)(14.380938,-0.19171876)
\psline[linewidth=0.04cm](14.380938,1.6082813)(16.380938,-0.19171876)
\psline[linewidth=0.072cm](2.740198,0.38611874)(3.0435576,0.42166913)
\psline[linewidth=0.072cm](2.740198,0.38611874)(2.7017639,0.6667184)
\psline[linewidth=0.072cm](7.407006,1.5908109)(7.2091036,1.7928872)
\psline[linewidth=0.072cm](7.407006,1.5908109)(7.2049294,1.3929089)
\psline[linewidth=0.072cm](1.4380041,-0.19841328)(1.658998,-0.3749401)
\psline[linewidth=0.072cm](1.4380041,-0.19841328)(1.6145309,0.022580562)
\psline[linewidth=0.072cm](1.8247576,0.77638924)(1.7555938,0.4788873)
\psline[linewidth=0.072cm](1.8247576,0.77638924)(2.1018164,0.71763813)
\psline[linewidth=0.072cm](15.198826,-0.17994644)(15.324124,-0.45849848)
\psline[linewidth=0.072cm](15.198826,-0.17994644)(15.454806,-0.058755502)
\psline[linewidth=0.072cm](0.7201978,0.38611874)(1.0235577,0.42166913)
\psline[linewidth=0.072cm](0.7201978,0.38611874)(0.6817638,0.6667184)
\psline[linewidth=0.072cm](8.382194,0.48254925)(8.633459,0.65620613)
\psline[linewidth=0.072cm](8.382194,0.48254925)(8.216663,0.71236)
\psline[linewidth=0.072cm](15.510726,0.6651254)(15.463106,0.94393075)
\psline[linewidth=0.072cm](15.510726,0.6651254)(15.231921,0.6175059)
\psline[linewidth=0.072cm](3.9358509,0.6037923)(4.008951,0.8770254)
\psline[linewidth=0.072cm](3.9358509,0.6037923)(3.6626177,0.67689264)
\psline[linewidth=0.072cm](2.3851204,1.5787519)(2.210053,1.8009037)
\psline[linewidth=0.072cm](2.3851204,1.5787519)(2.1629686,1.4036845)
\psline[linewidth=0.072cm](13.40512,1.5387518)(13.230053,1.7609036)
\psline[linewidth=0.072cm](13.40512,1.5387518)(13.182968,1.3636844)
\psline[linewidth=0.072cm](9.38512,1.5787519)(9.2100525,1.8009037)
\psline[linewidth=0.072cm](9.38512,1.5787519)(9.162969,1.4036845)
\psline[linewidth=0.072cm](9.380613,-0.2427602)(9.605781,-0.41393232)
\psline[linewidth=0.072cm](9.380613,-0.2427602)(9.551785,-0.017593354)
\psline[linewidth=0.072cm](7.380614,-0.2427602)(7.6057806,-0.41393232)
\psline[linewidth=0.072cm](7.380614,-0.2427602)(7.551786,-0.017593354)
\psline[linewidth=0.072cm](3.5780041,-0.15841329)(3.7989979,-0.3349401)
\psline[linewidth=0.072cm](3.5780041,-0.15841329)(3.754531,0.06258056)
\psline[linewidth=0.072cm](6.3941946,0.7393645)(6.1961656,0.50682336)
\psline[linewidth=0.072cm](6.3941946,0.7393645)(6.6133556,0.5599711)
\psline[linewidth=0.072cm](10.414195,0.77936447)(10.216166,0.5468233)
\psline[linewidth=0.072cm](10.414195,0.77936447)(10.633355,0.5999711)
\psline[linewidth=0.072cm](12.414195,0.7593645)(12.216166,0.52682334)
\psline[linewidth=0.072cm](12.414195,0.7593645)(12.633355,0.5799711)
\psline[linewidth=0.072cm](14.402193,0.46254927)(14.653459,0.6362061)
\psline[linewidth=0.072cm](14.402193,0.46254927)(14.236663,0.69236)
\usefont{T1}{ptm}{m}{n}
\rput(1.2323438,2.0182812){$3$}
\usefont{T1}{ptm}{m}{n}
\rput(3.3723438,1.9982812){$3$}
\usefont{T1}{ptm}{m}{n}
\rput(2.6123438,-0.68171877){$\geq3$}
\usefont{T1}{ptm}{m}{n}
\rput(4.532344,-0.70171875){$2$}
\usefont{T1}{ptm}{m}{n}
\rput(0.23234375,-0.74171877){$1$}
\usefont{T1}{ptm}{m}{n}
\rput(2.3782814,-1.8967187){\large $G_1$}
\usefont{T1}{ptm}{m}{n}
\rput(8.558281,-1.8767188){\large $G_2$}
\usefont{T1}{ptm}{m}{n}
\rput(14.258282,-1.9367187){\large $G_3$}
\psline[linewidth=0.072cm](13.2206135,-0.2227602)(13.445781,-0.3939323)
\psline[linewidth=0.072cm](13.2206135,-0.2227602)(13.391786,0.0024066463)
\usefont{T1}{ptm}{m}{n}
\rput(6.3323436,-0.70171875){$2$}
\usefont{T1}{ptm}{m}{n}
\rput(8.332344,-0.7217187){$\geq2$}
\usefont{T1}{ptm}{m}{n}
\rput(10.332344,-0.70171875){$2$}
\usefont{T1}{ptm}{m}{n}
\rput(10.332344,1.9982812){$2$}
\usefont{T1}{ptm}{m}{n}
\rput(6.2923436,1.9982812){$2$}
\usefont{T1}{ptm}{m}{n}
\rput(8.372344,1.9982812){$3$}
\usefont{T1}{ptm}{m}{n}
\rput(14.332344,-0.6617187){$\geq2$}
\usefont{T1}{ptm}{m}{n}
\rput(12.352344,-0.7217187){$2$}
\usefont{T1}{ptm}{m}{n}
\rput(16.472343,-0.6617187){$2$}
\usefont{T1}{ptm}{m}{n}
\rput(14.412344,2.0182812){$3$}
\usefont{T1}{ptm}{m}{n}
\rput(12.232344,1.9782813){$2$}
\end{pspicture} 
}
\caption{An orientation of the configurations }
\end{figure}
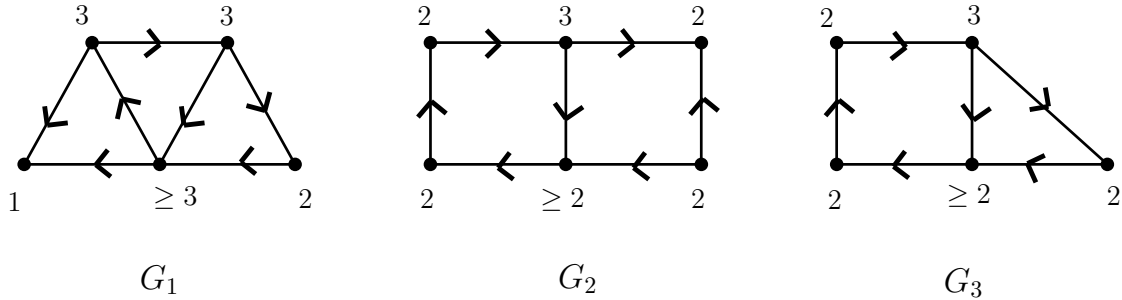

\begin{pr}
Let $L$ be a 4-list assignment of $G.$ 
Suppose $G$ contains a subgraph $H$ isomorphic to one of the configurations 
in Fig. \ref{fig3}. 
By minimality of $G,$ the graph $G -V(H)$ has an $L$-coloring. 
An orientation and availability of colors for every vertex 
are shown in Fig. \ref{fig4}. 
Since $EE(G_1)= 2 > EO(G_1)=1, EE(G_2)= 3 > EO(G_2)=1,$ and 
$EE(G_3)= 2 > EO(G_3)=1,$  each $G_i$ satisfies Theorem \ref{null}. 
Thus we can extend an $L$-coloring to $G,$ a contradiction. 
\end{pr}

\section{Main results and proof}

\begin{thm}\label{main1} 
If each 5-cycle of a planar graph $G$ is not a subgraph of 5-wheel 
and does not share exactly one edge with 3-cycle, 
then $G$ is 4-choosable. 
\end{thm}

\begin{thm}\label{main2} 
If each 5-cycle of a planar graph $G$ is not adjacent to two 3-cycles 
and is not adjacent to $\theta(1,2,2),$ then $G$ is 4-choosable. 
\end{thm} 

The following corollary is an easy consequence. 
\begin{cor}
If each 5-cycle of a planar graph $G$ is not adjacent a 3-cycle, 
then $G$ is 4-choosable. 
\end{cor} 

We use the same proof for Theorems \ref{main1} and \ref{main2}. 
In fact the proof of the latter is easier because it does not involve a trio. 

\begin{pr}

\indent Suppose that $G$ is a minimal counterexample. 
The discharging process is as follows. 
Let the initial charge of a vertex $v$ in $G$ be $\mu(v)=2d(v)-6$ 
and the initial charge of a face $f$ in $G$ be $\mu(f)=d(f)-6$. 
Then by Euler's formula $|V(G)|-|E(G)|+|F(G)|=2$ 
and by the Handshaking lemma, we have
$$\displaystyle\sum_{v\in V(G)}\mu(v)+\displaystyle\sum_{f\in F(G)}\mu(f)=-12.$$

\indent Now, we establish a new charge $\mu^*(x)$ 
for all $x\in V(G)\cup F(G)$ by transferring charge from one element to another
and the summation of new charge $\mu^*(x)$ remains $-12$. 
If the final charge  $\mu^*(x)\geq 0$ for all $x\in V(G)\cup F(G)$, 
then we get a contradiction and the proof is completed.\\

\indent Let $w(v \rightarrow f)$ be the charge transfered from 
a vertex $v$ to a face $f.$ 
The discharging rules are as follows.\\ 

\textbf{(R1)} $w(v \rightarrow f) = 0.2$ if $f$ is a 5-face.
 
\textbf{(R2)} For a 4-vertex $v$,  $w(v \rightarrow f) = 1$ if $v$ is 
 a good, bad, or, worse vertex of a 3-face $f,$ 
 $w(v \rightarrow f) = 2/3$ if $v$ is a worst vertex of a 3-face $f,$ 
 $w(v \rightarrow f) = 2/3$ if $v$ is a worst vertex of a 3-face $f,$ 
 and  $w(v \rightarrow f) = 1/3$ if $v$ is vertex on a 4-face $f.$

\textbf{(R3)} For a 5-vertex $v$,  $w(v \rightarrow f) = 1$ if $v$ is 
a good or worst vertex of a 3-face $f,$ 
$w(v \rightarrow f) = 1.5$ if $v$ is a bad vertex of a 3-face $f,$ 
$w(v \rightarrow f) = 1.25$ if $v$ is a worse vertex of a 3-face $f,$  
$w(v \rightarrow f) = 1$ if $v$ is vertex on a (4,4,4,5)-face $f,$ 
$w(v \rightarrow f) = 2/3$ if $v$ is vertex on a 4-face $f$ 
that is not a (4,4,4,5)-face.

\textbf{(R4)} For a 6$^+$-vertex $v$,  $w(v \rightarrow f) = 1$ if $v$ is 
a good or worst vertex of a 3-face $f,$ 
$w(v \rightarrow f) = 1.5$ if $v$ is a bad vertex of a 3-face $f,$ 
$w(v \rightarrow f) = 1.25$ if $v$ is a worse vertex of a 3-face $f,$  
$w(v \rightarrow f) = 1$ if $v$ is vertex on a (4,4,4,5)-face $f,$ 
$w(v \rightarrow f) = 2/3$ if $v$ is vertex on a 4-face $f$ 
that is not a (4,4,4,5)-face. 

\textbf{(R5)} After (R1)-(R4), redistribute charge in three 3-faces 
in the same trio to make their charge equal.\\ 

\indent It remains to show that resulting $\mu^*(x)\geq 0$ 
for all $x\in V(G)\cup F(G)$.\\

Consider a 4-vertex $v.$ If $v$ is not incident to any 3-face, 
then $\mu^*(v)\geq \mu(v) - 4\cdot (1/3) = 2/3.$ 
If $v$ is  incident to exactly one 3-face, 
then  $\mu^*(v)\geq \mu(v) - (1 + 3\cdot (1/3)) = 0.$ 
If $v$ is incident to  two 3-faces, then $v$ is not a worst vertex of any face 
and  other incident faces are $6^+.$ 
Using (R2), we obtain  $\mu^*(v) = \mu(v) - (2\cdot 1)= 0.$ 
If $v$ is incident to  three 3-faces, 
then $v$ is a worst vertex of these faces and the other incident face is $6^+.$ 
Using (R2), we obtain  $\mu^*(v) = \mu(v) - (3\cdot (2/3)) = 0.$ 

Consider a 5-vertex $v.$ 
If $v$ is a (3,4,4,4,5)- or a (4,4,4,4,5)-vertex, 
then at least one of its incident faces is a (4,4,5,5)-face by Lemma \ref{lem2}(3). 
This means $v$ sends charge to its incident faces at most 
$3+(2/3)+(1/5)$ which is less than 4. 
If $v$ is a (3,4,4,4,4)- or a (4,4,4,4,4)-vertex, 
then at least three of its incident faces is a (4,4,5,5)-face by Lemma \ref{lem2}(3). 
This means $v$ sends charge to its incident faces at most 
$2+3\cdot(2/3)=4.$ 
If $v$ is incident to at least two 5-face, 
then $v$ sends charge at most $2\cdot(1/5)+3$ in the case that 
$v$ has a bad face, 
and at most $2\cdot(1/5)+3$ if $v$ has no bad faces. 
If $v$ has two bad faces, then $v$ has three 6$^+$-faces.  
Using (R3), we obtain  $\mu^*(v) = \mu(v) - (2\cdot (1.5)) = 2.$ 
So we assume that $v$ is incident to a 6$^+$-face and has at most one bad face. 
If $v$ has a worst face, a worse face, or a bad face, then 
$v$ is incident to at least two 6$^+$-faces 
and $\mu^*(v) = \mu(v) - 3 = 1$ in the case that $v$ has a worst face,  
or $\mu^*(v) \geq \mu(v) - (2\cdot (1.25)+1) >0$ in the case that 
$v$ has  a worse face or a bad face. 
If $v$ has neither a worst face, a worse face, nor a bad face, 
then $\mu^*(v) \geq \mu(v) - 4 =0.$ 

If $v$ has a worst face, a worse face or a bad face, respectively, then 
$v$ is incident to at least two 6$^+$-faces 
and $\mu^*(v) = \mu(v) - 3 = 1$ 
or $\mu^*(v) \geq \mu(v) - (2\cdot (1.25)+1) >0,$ respectively.

Consider a $k$-vertex $v$ with $k\geq 6.$ 
Assume $v$ is incident to the faces $f_1,\ldots,f_k$ in a cyclic order,  
where all subscripts are taken modulo $k.$ 
To calculate $\mu(v),$  we redistribute $w(v \rightarrow f_i).$ 
However, we still use  $w(v \rightarrow f_i)$ according to (R4) 
for the final charge.  
If $w(v \rightarrow f_i) = 1.5,$  
then $v$ is a worse vertex of $f_i$ and $f_{i-1}$ (or $f_{i+1}$) is a 6$^+$-face. 
Reduce $w(v \rightarrow f_i)$ to 1 and 
transfer the remaining charge 0.5 to $f_{i-1}.$ 
If $w(v \rightarrow f_i) = 2,$  
then $v$ is a bad vertex of $f_i$ whereas both $f_{i-1}$ and $f_{i+1}$) 
are 6$^+$-faces. 
Reduce $w(v \rightarrow f_i)$ to 1 and transfer charge 0.5 to $f_{i-1}$ 
and $f_{i+1}.$ According to this average, 
$v$ sends charge at most 1 to each of its adjacent faces. 
Thus $\mu^*(v) \geq \mu(v) - \deg(v) \geq 0.$

\indent It is clear  that $\mu^*(f)=\mu(f) \geq 0$ for $f$ a $6^+$-face $f.$  
If $f$ is a 5-face, then $\mu^*(f)=\mu(f) +5\cdot(0.2) =0$ by (R1). 
Consider a 4-face $f.$ 
Lemma \ref{lem1} yields that every vertex on $f$ has degree at least 4 and one 
of them has degree at least 5. 
If $f$ is a $(4,4,4,5)$-face, 
then $\mu^*(f)\geq \mu(f) + 3(2/3)+1 =0,$   
otherwise  $\mu^*(f) \geq \mu(f) +2(1/3)+2 \cdot (2/3) =1$ by (R2) and (R3).  
If $f$ is a 3-face that is not in a trio, then  $\mu^*(f) = \mu(f) +3(1/3) =0$ 
by (R2), (R3), and (R4). 

According to (R5), it remains to show that the total increased charge 
of three 3-faces in each trio is at least 9 
(or in the other words, the summation of their new charge is 0.)   
If the degree of a worst vertex is at least 5, 
then each vertex send charge at least 1 to its incident face. 
So we assume that a worst vertex has degree 4 and it sends charge totaling 
$3\cdot (2/3) =2$ to three faces in a trio. 
By (R2)-(R4), each of the remaining vertices in the trio sends charge at least 1 
to its incident face in a trio. 
Note that  a 6$^+$-vertex $v$ in the trio sends charge totaling $2 \cdot (1.5)= 3$ 
to two faces in the trio if $v$ is a worse vertex, 
otherwise it sends charge 2 to a face in the trio according to (R4). 
Additionally, a 5-vertex $v$ in a trio sends charge totaling $2 \cdot (1.25)= 2.5$ 
to two faces in the trio if $v$ is a worse vertex, 
otherwise it sends charge 1.5 to a face in the trio according to (R4). 
Thus total increased charge is at least 9 if the trio has a  6$^+$-vertex $v$ 
or two $5$-vertex. 
Use Lemma \ref{lem1} and \ref{lem2} to eliminate the remaining case. 
This completes the proof. 
\end{pr}

It is interesting to see that whether one can reduce some restrictions in 
Theorems \ref{main1} or \ref{main2} but a planar graph is still 4-choosable.

\section*{Acknowledgments} 
The first author is supported by Development and Promotion of Science 
and Technology talents project (DPST).



\begin{thebibliography}{99}
\bibitem{Alon}
N. Alon, M. Tarsi, Colorings and orientations of graphs, Combinatorica 12(1992) 125-134.

\bibitem{Bo}
O.V. Borodin, A.O. Ivanova, Planar graphs without triangular 4-cycles are 4-choosable, Sib. \' Elektron. Mat. Rep. 5(2008) 75-79.(check jounal 's name

\bibitem{Cheng}
P.P. Cheng, M. Chen, Y.Q. Wang, 
Planar graphs without 4-cycles adjacent to triangles are 4-choosable, 
Discrete Math. 339(2016) 3052-3057.\\

\bibitem{Erdos}
P. Erd\H os, A.L. Rubin, H. Taylor, 
Choosability in graphs, in: 
Proceedings, West Coast Conference on Combinatorics, 
Graph Theory and Computing, Arcata,
CA., Sept. 5-7, in: Congr. Numer., vol. 26, 1979.

\bibitem{Far}
B. Farzad, Planar graphs without 7-cycles are 4-choosable, 
SIAM J. Discrete Math. 23(2009) 1179-1199.

\bibitem{Fi}
G. Fijav\v z, M. Juvan, B. Mohar, R. \v Skrekovski, 
Planar graphs without cycles of specific lengths, European J. Combin. 23(2002) 377-388.

\bibitem{Gut}
S. Gutner, The complexity of planar graph choosability, 
Discrete Math. 159(1996) 119-130.

\bibitem{Hu}
D.Q. Hu, J.L. Wu, Planar graphs without intersecting 5-cycles are 4-choosable, 
Discrete Math. 340(2017) 1788-1792. 

\bibitem{Lam1}
P.C.B. Lam, W.C. Shiu, B.G. Xu, On structure of some plane graphs with applications to choosability, J. Combin. Theory Ser. B 82(2001) 285-296.

\bibitem{Lam2} 
P.C.B. Lam, B. Xu, J. Liu, The 4-choosability of plane graphs without 4-cycles, J. Combin. Theory Ser. B 76(1999) 117-126.

\bibitem{Mir}
M. Mirzakhani, A small non-4-choosable planar graph, Bull. Inst. Combin. Appl. 17(1996) 15-18.

\bibitem{Tho}
C. Thomassen, Every planar graph is 5-choosable, 
J. Combin. Theory Ser. B 62(1994) 180-181.

\bibitem{Vizing} 
V.G. Vizing, Vertex colorings with given colors, 
Metody Diskret. Anal. 29(1976) 3-10 (in Russian).

\bibitem{Vo1}
M. Voigt, List colourings of planar graphs, Discrete Math. 120(1993) 215-219.


\bibitem{Xu} 
R. Xu, J.L. Wu, 
A sufficient condition for a planar graph to be 4-choosable, 
Discrete App. Math. 224(2017)120-122.

\bibitem{Wang1}
W. Wang, K.W. Lih, 
Choosability and edge choosability of planar graphs without five cycles, 
Appl. Math. Lett. 15(2002) 561-565.

\bibitem{Wang2} 
W. Wang, K.W. Lih, 
Choosability and edge choosability of planar without intersecting triangles, 
SIAM J. Discrete Math. 15(2002) 538-545.

\end{thebibliography}
\end{document}